\documentstyle[12pt]{article}
\begin{document}
\begin{titlepage}
\begin{flushright}
math.QA/9811052
\end{flushright}
\vskip.3in

\begin{center}
{\Large \bf Casimir Invariants from Quasi-Hopf (Super)algebras}
\vskip.3in
{\large Mark. D. Gould, Yao-Zhong Zhang} and {\large Phillip S. Isaac}
\vskip.2in
{\em Department of Mathematics, University of Queensland, Brisbane,
     Qld 4072, Australia

Email: yzz@maths.uq.edu.au}
\end{center}

\vskip 2cm
\begin{center}
{\bf Abstract}
\end{center}

We show how to construct, starting from a quasi-Hopf (super)algebra,
central elements or Casimir invariants. We show that these central
elements are invariant under quasi-Hopf twistings. As a consequence,
the elliptic quantum (super)groups, which arise from twisting the
normal quantum (super)groups, have the same Casimir invariants as
the corresponding quantum (super)groups.

\vskip 3cm
\noindent{\bf Mathematics Subject Classifications (1991):} 81R10, 17B37, 16W30

\end{titlepage}


\def\a{\alpha}
\def\b{\beta}
\def\d{\delta}
\def\e{\epsilon}
\def\ve{\varepsilon}
\def\g{\gamma}
\def\k{\kappa}
\def\l{\lambda}
\def\o{\omega}
\def\t{\theta}
\def\s{\sigma}
\def\D{\Delta}
\def\L{\Lambda}

\def\G{{\cal G}}
\def\hG{{\hat{\cal G}}}
\def\R{{\cal R}}
\def\hR{{\hat{\cal R}}}
\def\C{{\bf C}}
\def\P{{\bf P}}
\def\Z2{{{\bf Z}_2}}
\def\T{{\cal T}}
\def\H{{\cal H}}
\def\trho{{\tilde{\rho}}}
\def\tphi{{\tilde{\phi}}}
\def\tT{{\tilde{\cal T}}}
\def\uqsnh{{U_q[\widehat{sl(n|n)}]}}
\def\uqs1h{{U_q[\widehat{sl(1|1)}]}}


\def\beq{\begin{equation}}
\def\eeq{\end{equation}}
\def\bea{\begin{eqnarray}}
\def\eea{\end{eqnarray}}
\def\ba{\begin{array}}
\def\ea{\end{array}}
\def\no{\nonumber}
\def\lt{\left}
\def\rt{\right}
\newcommand{\bq}{\begin{quote}}
\newcommand{\eq}{\end{quote}}

\newtheorem{Theorem}{Theorem}
\newtheorem{Definition}{Definition}
\newtheorem{Proposition}{Proposition}
\newtheorem{Lemma}{Lemma}
\newtheorem{Corollary}{Corollary}
\newcommand{\proof}[1]{{\bf Proof. }
        #1\begin{flushright}$\Box$\end{flushright}}

\newcommand{\sect}[1]{\setcounter{equation}{0}\section{#1}}
\renewcommand{\theequation}{\thesection.\arabic{equation}}

\sect{Introduction\label{intro}}

Quasi-Hopf superalgebras are $\Z2$-graded versions of Drinfeld's
quasi-Hopf algebras \cite{Dri90} and were introduced in \cite{Zha98}.
The potential for applications of these structures, particularly to knot
theory and integrable systems, is enormous. They give rise to new
(non-standard) representations of the braid group and corresponding
link polynomials \cite{Alt92,Mac91}. Moreover these remarkable structures
underly elliptic quantum (super)groups 
\cite{Fod94,Fel95,Fro97,Jim97,Enr97,Zha98} which play an important
role in obtaining solutions to the dynamical Yang-Baxter equations
\cite{Bab96,Arn97}.

In applications such as these it is important to have a well defined
representation theory. In this paper we investigate several aspects
of this theory concerned with the construction and general properties
of invariants (invariant bilinear forms, module morphisms, central
elements and {\it etc}). In particular, in the quasi-triangular case,
it is shown how central elements may be systematically constructed
utilizing the R-matrix. This construction may be regarded as a
natural generalization of that introduced in \cite{Lin92,Zha91}, 
to which it reduces
in the case of normal Hopf (super)algebras. However the extension of
this paper is by no means straightforward and requires the explicit
inclusion of the co-associator into the construction.

We moreover prove the strong result that the Casimir invariants so
obtained are invariant under twisting. This implies, in particular,
that one will not obtain new Casimir invariants by twisting on
quantum (super)groups. As part of our approach we extend the
$u$-operator formalism of Drinfeld--Reshetikhin to the case of
quasi-Hopf superalgebras. In particular we prove the surprising
result that the $u$-operator is invariant under twisting. This has some
important implications for knot theory which will be investigated
elsewhere. It is worth noting that most of our results are new, even
in the non-graded case.

\sect{Quasi-Hopf (Super)algebras}

Let us briefly recall the quasi-Hopf algebras \cite{Dri90} and their
super (or $\Z2$ graded) versions -- quasi-Hopf superalgebras
\cite{Zha98}.

\begin{Definition}\label{quasi-bi}: 
A quasi-Hopf (super)algebra is a ($\Z2$ graded) 
unital associative algebra $A$
over a field $K$  which is equipped with algebra homomorphisms $\e: 
A\rightarrow K$ (co-unit), $\D: A\rightarrow A\otimes A$ (co-product),
an invertible homogeneous element $\Phi\in A\otimes A\otimes A$ 
(co-associator), an 
($\Z2$ graded) algebra anti-homomorphism $S: A\rightarrow A$ (anti-pode) and 
homogeneous canonical elements $\a,~\b\in A$, satisfying 
\bea
&& (1\otimes\D)\D(a)=\Phi^{-1}(\D\otimes 1)\D(a)\Phi,~~
	\forall a\in A,\label{quasi-bi1}\\
&&(\D\otimes 1\otimes 1)\Phi \cdot (1\otimes 1\otimes\D)\Phi
	=(\Phi\otimes 1)\cdot(1\otimes\D\otimes 1)\Phi\cdot (1\otimes
	\Phi),\label{quasi-bi2}\\
&&(\e\otimes 1)\D=1=(1\otimes\e)\D,\label{quasi-bi3}\\
&&(1\otimes\e\otimes 1)\Phi=1,\label{quasi-bi4}\\
&& m\cdot (1\otimes\a)(S\otimes 1)\D(a)=\e(a)\a,~~~\forall
    a\in A,\label{quasi-hopf1}\\
&& m\cdot (1\otimes\b)(1\otimes S)\D(a)=\e(a)\b,~~~\forall a\in A,
     \label{quasi-hopf2}\\
&& m\cdot (m\otimes 1)\cdot (1\otimes\b\otimes\a)(1\otimes S\otimes
     1)\Phi^{-1}=1,\label{quasi-hopf3}\\
&& m\cdot(m\otimes 1)\cdot (S\otimes 1\otimes 1)(1\otimes\a\otimes
     \b)(1\otimes 1\otimes S)\Phi=1.\label{quasi-hopf4}
\eea
\end{Definition}
Here $m$ denotes the usual product map on $A$: $m\cdot (a\otimes b)=ab,~
\forall a,b\in A$. Note that since $A$ is associative we have
$m\cdot(m\otimes 1)=m\cdot (1\otimes m)$.
For the homogeneous elements $a,b\in A$, the antipode satisfies
\beq
S(ab)=(-1)^{[a][b]}S(b)S(a),
\eeq
which extends to inhomogeneous elements through linearity.
(\ref{quasi-bi2}), (\ref{quasi-bi3}) and (\ref{quasi-bi4}) imply
that $\Phi$ also obeys
\beq
(\e\otimes 1\otimes 1)\Phi=1=(1\otimes 1\otimes\e)\Phi.\label{e(phi)=1}
\eeq
It follows that the co-associator $\Phi$ is an even element.
Applying $\e$ to definition (\ref{quasi-hopf3}, \ref{quasi-hopf4}) 
we obtain, in view
of (\ref{quasi-bi4}), $\e(\a)\e(\b)=1$. Thus the canonical
elements $\a, \b$ are both even. By applying $\e$ to
(\ref{quasi-hopf1}), we have $\e(S(a))=\e(a),~\forall a\in A$.
Note that the multiplication rule for the tensor products is defined
for homogeneous elements $a,b,a',b'\in A$ by
\beq
(a\otimes b)(a'\otimes b')=(-1)^{[b][a']}\,(aa'\otimes bb'),
\eeq
where $[a]\in{\bf Z}_2$ 
denotes the grading of the element $a$.

The category of quasi-Hopf (super)algebras
is invariant under a kind of gauge transformation. Let $(A,\D,\e,\Phi)$
be a qausi-Hopf (super)algebra, with $\a,\b, S$ satisfying
(\ref{quasi-hopf1})-(\ref{quasi-hopf4}), and let $F\in A\otimes A$
be an invertible homogeneous element satisfying the co-unit properties
\beq
(\e\otimes 1)F=1=(1\otimes \e)F.\label{e(f)=1}
\eeq
It follows that $F$ is even. Throughout we set
\bea
&&\D_F(a)=F\D(a)F^{-1},~~~\forall a\in A,\label{twisted-d}\\
&&\Phi_F=(F\otimes 1)(\D\otimes
    1)F\cdot\Phi\cdot(1\otimes\D)F^{-1}(1\otimes F^{-1}).\label{twisted-phi}
\eea 
Then 
\begin{Theorem}\label{t-quasi-hopf}:
$(A,\D_F,\e,\Phi_F)$ defined by (\ref{twisted-d}, 
\ref{twisted-phi}) together with
$\a_F,\b_F, S_F$ given by
\beq
S_F=S,~~~\a_F=m\cdot(1\otimes\a)(S\otimes 1)F^{-1},~~~
	 \b_F=m\cdot(1\otimes\b)(1\otimes S)F,\label{twisted-s-ab}
\eeq
is also a quasi-Hopf (super)algebra. The element $F$ is referred to as
a twistor, throughout.
\end{Theorem}

\begin{Definition}\label{quasi-quasi}: A quasi-Hopf
(super)algebra $(A,\D,\e,\Phi)$ is called quasi-triangular if there
exists an invertible homogeneous element $\R\in A\otimes A$ such that
\bea
&&\D^T(a)\R=\R\D(a),~~~~\forall a\in A,\label{dr=rd}\\
&&(\D\otimes 1)\R=\Phi^{-1}_{231}\R_{13}\Phi_{132}\R_{23}\Phi^{-1}_{123},
   \label{d1r}\\
&&(1\otimes \D)\R=\Phi_{312}\R_{13}\Phi^{-1}_{213}\R_{12}\Phi_{123}.
   \label{1dr}
\eea
$\R$ is referred to as the universal R-matrix.
\end{Definition}
Throughout, $\D^T=T\cdot\D$ with $T$ being the graded twist map
which is defined, for homogeneous elements $a,b\in A$, by
\beq
T(a\otimes b)=(-1)^{[a][b]}b\otimes a;
\eeq
and $\Phi_{132}$ {\it etc} are derived from $\Phi\equiv\Phi_{123}$
with the help of $T$
\bea
&&\Phi_{132}=(1\otimes T)\Phi_{123},\no\\
&&\Phi_{312}=(T\otimes 1)\Phi_{132}=(T\otimes 1)
   (1\otimes T)\Phi_{123},\no\\
&&\Phi^{-1}_{231}=(1\otimes T)\Phi^{-1}_{213}=(1\otimes T)
   (T\otimes 1)\Phi^{-1}_{123},\no
\eea
and so on. 

It is easily shown that the properties (\ref{dr=rd})-(\ref{1dr})
imply the (graded) Yang-Baxter type equation,
\beq
\R_{12}\Phi^{-1}_{231}\R_{13}\Phi_{132}\R_{23}\Phi^{-1}_{123}
  =\Phi^{-1}_{321}\R_{23}\Phi_{312}\R_{13}\Phi^{-1}_{213}\R_{12},
  \label{quasi-ybe}
\eeq
which is referred to as the (graded)  quasi-Yang-Baxter equation, 
and the co-unit properties of $\R$:
\beq
(\e\otimes 1)\R=1=(1\otimes \e)\R.\label{e(R)=1}
\eeq
Thus the universal R-matrix $\R$ is even. We have

\begin{Theorem}\label{t-quasi-quasi}: Denoting by the set
$(A,\D,\e,\Phi,\R)$  a
quasi-triangular quasi-Hopf (super)algebra, then $(A, \D_F, \e, \Phi_F, \R_F)$
is also a quasi-triangular quasi-Hopf (super)algebra, with the choice of
$R_F$ given by
\beq
\R_F=F^T \R F^{-1},\label{twisted-R}
\eeq
where $F^T=T\cdot F\equiv F_{21}$. Here $\D_F$ and $\Phi_F$ are given
by (\ref{twisted-d}) and (\ref{twisted-phi}), respectively.
\end{Theorem}

Let us specify some notations. Throughtout the paper,
\bea
&&\Phi=\sum X_\nu\otimes Y_\nu\otimes Z_\nu,~~~~
  \Phi^{-1}=\sum \bar{X}_\nu\otimes\bar{Y}_\nu\otimes\bar{Z}_\nu,\no\\
&&F=\sum f_i\otimes f^i,~~~~F^{-1}=\sum \bar{f}_i\otimes \bar{f}^i,\no\\
&&\R=\sum e_i\otimes e^i,~~~~\R^{-1}=\sum\bar{e_i}\otimes\bar{e}^i,\no\\
&&(1\otimes\D)\D(a)=\sum a_{(1)}\otimes\D(a_{(2)})=\sum
   a^R_{(1)}\otimes a_{(2)}^R\otimes a^R_{(3)},\no\\
&&(\D\otimes 1)\D(a)=\sum \D(a_{(1)})\otimes a_{(2)})=\sum
   a^L_{(1)}\otimes a_{(2)}^L\otimes a^L_{(3)}.\label{notation}
\eea
The following lemma is proved in \cite{Gou98} and will be used
frequently in this paper.
\begin{Lemma}\label{L1}: $\forall a\in A$,
\bea
&(i)& \sum X_\nu a\otimes Y_\nu\b S(Z_\nu)(-1)^{[a][X_\nu]}
      =\sum a^L_{(1)}X_\nu\otimes a^L_{(2)}Y_\nu\b S(Z_\nu)S(a^L_{(3)})
      (-1)^{[a^L_{(2)}][X_\nu]},\no\\
&(ii)& \sum S(X_\nu)\a Y_\nu\otimes a Z_\nu (-1)^{[a][Z_\nu]}=
     \sum S(a^R_{(1)})S(X_\nu)\a Y_\nu a^R_{(2)}\otimes Z_\nu a^R_{(3)}
     (-1)^{[a^R_{(2)}][Z_\nu]},\no\\
&(iii)& \sum a \bar{X}_\nu \otimes S(\bar{Y}_\nu)\a \bar{Z}_\nu
      =\sum \bar{X}_\nu a^L_{(1)}\otimes S(a^L_{(2)})S(\bar{Y}_\nu)\a 
      \bar{Z}_\nu a^L_{(3)}(-1)^{[X_\nu]([a^L_{(1)}]+[a^L_{(2)}])},\no\\
&(iv)& \sum \bar{X}_\nu\b S(\bar{Y}_\nu)\otimes \bar{Z}_\nu a=
     \sum a^R_{(1)}\bar{X}_\nu\b S(\bar{Y}_\nu)S(a^R_{(2)})\otimes a^R_{(3)}
     \bar{Z_\nu}(-1)^{([a^R_{(2)}]+[a^R_{(3)}])[\bar{Z}_\nu]}.
\eea
\end{Lemma}

\sect{Central Elements from (Anti-)adjoint Actions}

Given an ($\Z2$ graded) $A$-module $V$, we say $v\in V$ an invariant if
\beq
a\cdot v=\e(a) v,~~~~\forall a\in A.
\eeq
In particular, $A$ itself constitutes an $A$-module under the adjoint
action defined by
\beq
{\rm Ad} a\cdot b=\sum a_{(1)} b S(a_{(2)}) (-1)^{[b][a_{(2)}]},~~~~
 \forall a,b\in A.
\eeq
It is easily shown that 
\beq
{\rm Ad}a\cdot{\rm Ad} b={\rm Ad} ab.
\eeq
We call $c_1\in A$ an invariant if it is invariant under the adjoint
action, i.e.
\beq
\sum a_{(1)} c_1S(a_{(2)}) (-1)^{[c_1][a_{(2)}]}=\e(a) c_1,~~~~
   \forall a\in A.\label{c-inv}
\eeq
For normal Hopf (super)algebras, the invariants of $A$ are precisely the
central elements. This is not true, however, for quasi-Hopf
(super)algebras. For instance, the canonical element $\b$ is
invariant but not generally central. Nevertheless, there is a close
connection between central elements and invariants. We have
\begin{Proposition}\label{P-cbar1}: 
Suppose $c_1\in A$ is even and invariant. Set
\beq
C_1=\sum \bar{X}_\nu c_1 S(\bar{Y}_\nu)\a\bar{Z}_\nu=m(m\otimes 1)\cdot
   (1\otimes c_1\otimes\a)(1\otimes S\otimes 1)\Phi^{-1}.\label{C1-1}
\eeq
Then $(i)~ a C_1=C_1a,~\forall a\in A$, i.e. $C_1$ is
central, and 
\bea
&(ii)& c_1=C_1\b=\b C_1,\no\\
&(iii)& C_1=\sum S(X_\nu)\a Y_\nu c_1 S(Z_\nu).\label{C1-2}
\eea
\end{Proposition}
\vskip.1in
\noindent{\bf Proof.}
Applying $m\cdot(1\otimes c_1)$ to Lemma \ref{L1}(ii) and
keeping in mind of (\ref{c-inv}), we obtain
(i). We now prove (ii). From (\ref{quasi-bi2}),
\bea
C_1\otimes 1&=&(m(m\otimes 1)\otimes 1)\cdot (1\otimes c_1\otimes \a
  \otimes 1)(1\otimes S\otimes 1\otimes 1) (\Phi^{-1}\otimes 1)\no\\
&=&\sum\lt(X_\nu\bar{X}_\s\bar{X}^{(1)}_\rho c_1S(\bar{X}^{(2)}_\rho)
  S(\bar{Y}_\s)S(X_\mu)S(Y^{(1)}_\nu)\a Y^{(2)}_\nu Y_\mu
  \bar{Z}^{(1)}_\s\bar{Y}_\rho\rt.\no\\
& &\lt.  \otimes Z_\nu Z_\mu\bar{Z}^{(2)}_\s
  \bar{Z}_\rho\rt) (-1)^x,
\eea
where
\bea
x&=&[\bar{Z}_\s][X_\nu]+[X_\mu][Y_\nu]+[Z_\mu][Z_\nu]+
  ([\bar{Z}^{(1)}_\s]+[\bar{Y}_\rho])([Z_\mu]+[Z_\nu])\no\\
& &  +[\bar{Y}_\rho][\bar{Z}^{(2)}_\s]+[\bar{X}_\rho]([X_\nu]+[\bar{X}_\s])
\eea
By (\ref{c-inv}) and (\ref{quasi-hopf1}),
\bea
C_1\otimes 1&=&
  \sum\lt(X_\nu\bar{X}_\s\e(\bar{X}_\rho) c_1
  S(\bar{Y}_\s)S(X_\mu)\e(Y_\nu)\a Y_\mu
  \bar{Z}^{(1)}_\s\bar{Y}_\rho\otimes Z_\nu Z_\mu\bar{Z}^{(2)}_\s
  \bar{Z}_\rho\rt) (-1)^x\no\\
&=&\sum (X_\nu\e(Y_\nu)\otimes Z_\nu)(
  \bar{X}_\s c_1
  S(\bar{Y}_\s)S(X_\mu)\a Y_\mu
  \bar{Z}^{(1)}_\s\otimes Z_\mu\bar{Z}^{(2)}_\s)\no\\
& &  (\e(\bar{X}_\rho)\bar{Y}_\rho\otimes\bar{Z}_\rho)(-1)^{[Z_\mu]
  [\bar{Z}^{(1)}_\s]}\no\\
&=&\sum(  \bar{X}_\s c_1
  S(\bar{Y}_\s)S(X_\mu)\a Y_\mu
  \bar{Z}^{(1)}_\s\otimes Z_\mu\bar{Z}^{(2)}_\s)(-1)^{[Z_\mu]
  [\bar{Z}^{(1)}_\s]}\no\\
& &  {\rm by}~(\ref{quasi-bi4}), (\ref{e(phi)=1}).
\eea
Applying $m\cdot(1\otimes\b)(1\otimes S)$ gives rise to
\bea
C_1\b&=&
  \sum \bar{X}_\s c_1
  S(\bar{Y}_\s)S(X_\mu)\a Y_\mu
  \bar{Z}^{(1)}_\s\b S(\bar{Z}^{(2)}_\s)S(Z_\mu)(-1)^{[Z_\mu]
  [\bar{Z}_\s]}\no\\
&=&  \sum \bar{X}_\s c_1
  S(\bar{Y}_\s)S(X_\mu)\a Y_\mu
  \e(\bar{Z}_\s)\b S(Z_\mu) ~~{\rm by}~(\ref{quasi-hopf2})\no\\
&=&  \sum \bar{X}_\s c_1
  S(\bar{Y}_\s)\e(\bar{Z}_\s)=c_1
  ,~~{\rm by}~(\ref{quasi-hopf4}), (\ref{e(phi)=1}),
\eea
thus proving (ii). (iii) is the direct consequence of (i) and (ii).

The above gives a very clear picture of the connection between
invariants and central elements. In particular we have
\begin{Corollary}: Suppose $c\in A$ is even. Then $c$ is an invariant
if and only if there exists a central element $C$ such that
\beq
c=C\b=\b C.
\eeq
\end{Corollary}

$A$ also admits an anti-adjoint action defined by 
\beq
\overline{\rm Ad}a\cdot b=\sum S(a_{(1)})ba_{(2)}(-1)^{[b][a_{(1)}]},~~~~
   \forall a, b\in A.
\eeq
We have
\beq
\overline{\rm Ad}a\cdot\overline{\rm Ad}b=\overline{\rm Ad}(ba).
\eeq
We call $c_2\in A$ a pseudo-invariant if it is invariant under the
anti-adjoint action; i.e.
\beq
\sum S(a_{(1)})c_2 a_{(2)}(-1)^{[c_2][a_{(1)}]}=\e(a)c_2,~~~~\forall
  a\in A.\label{c-anti-inv}
\eeq
\begin{Proposition}\label{P-cbar2}: 
Suppose $c_2\in A$ ie even and pseudo-invariant.
Set
\beq
C_2=\sum S(X_\nu)c_2 Y_\nu\b S(Z_\nu)=m(m\otimes 1)\cdot
   (1\otimes c_2\otimes\b)(S\otimes 1\otimes S)\Phi.\label{C2-1}
\eeq
Then $(i)~ aC_2=C_2a,~\forall a\in A$, i.e. $C_2$ is central, and
\bea
&(ii)& c_2=C_2\a=\a C_2,\no\\
&(iii)& C_2=\sum \bar{X}_\nu\b S(\bar{Y}_\nu) c_2 \bar{Z}_\nu.\label{C2-2}
\eea
\end{Proposition}
\vskip.1in
\noindent{\bf Proof.}
Similar to the proof of proposition \ref{P-cbar1}. 
Applying $m\cdot(1\otimes c_2)
(S\otimes 1) $ to Lemma \ref{L1}(i), we obtain (i). Applying 
$m\cdot(1\otimes \a)$ to 
\beq
C_2\otimes 1=(m(m\otimes 1)\otimes 1)\cdot (1\otimes c_2\otimes\b)
  (S\otimes 1\otimes S)(\Phi\otimes 1)
\eeq
leads to (ii). Finally, (iii) is a direct consequence of (i) and (ii).

As an example we construct the so-called quadratic invariants.
Suppose $\o=\sum\o_i\otimes\o^i\in A\otimes A$ is even and satisfies
\beq
\D(a)\o=\o\D(a),~~~~~\forall a\in A
\eeq
Applying $m\cdot(1\otimes\b)(1\otimes S)$ gives
\beq
\sum a_{(1)}\o_i\b S(\o^i) S(a_{(2)})=\e(a)\sum\o_i\b S(\o^i),~~~~
   \forall a\in A,
\eeq
which implies that
\beq
c_1\equiv \sum\o_i\b S(\o^i)\label{quadra-c1}
\eeq
is an invariant. Similarly, applying $m\cdot(1\otimes\a)(S\otimes 1)$,
one can show that
\beq
c_2\equiv \sum S(\o_i)\a\o^i\label{quadra-c2}
\eeq
is a pseudo-invariant. It follows from propositions \ref{P-cbar1} and
\ref{P-cbar2} that
\bea
C_1&=&\sum\bar{X}_\nu\o_i\b S(\o^i)S(\bar{Y}_\nu)\a\bar{Z}_\nu=
   \sum S(X_\nu)\a Y_\nu\o_i\b S(\o^i) Z_\nu,\no\\
C_2&=&\sum S(X_\nu)S(\o_i)\a\o^i Y_\nu\b S(Z_\nu)=\sum
   \bar{X}_\nu\b S(\bar{Y}_\nu)
   S(\o_i)\a\o^i\bar{Z}_\nu\label{quadra-cbars1}
\eea
are central elements. The invariants (\ref{quadra-cbars1}) are 
usually referred to as quadratic invariants.

\sect{Twisting Invariance of Central Elements $C_1$ and $C_2$}

\begin{Lemma}\label{L2}: Let $c_1\in A$ be even and invariant, and
$c_2\in A$ be even and pseudo-invariant. For any $\eta\in A\otimes A$,
$\xi\in A\otimes A\otimes A$, we have, $\forall a,b\in A$,
\bea
&(i)& m\cdot(1\otimes c_1)(1\otimes S)(\eta\D(a))=\e(a)
   m\cdot (1\otimes c_1)(1\otimes S)\eta,\no\\
&(ii)& m\cdot(1\otimes c_2)(S\otimes 1)(\D(a)\eta)=\e(a)
   m\cdot (1\otimes c_2)(S\otimes 1)\eta,\no\\
&(iii)& m(m\otimes c_1\otimes c_2) (1\otimes S\otimes 1)[(1\otimes
  \D(a))\cdot\xi\cdot(\D(b)\otimes 1)]\no\\
& &~~~~~~~~~~=\e(a)\e(b)m(m\otimes 1)\cdot (1\otimes c_1\otimes c_2)
  (1\otimes S\otimes 1)\xi,\no\\
&(iv)& m(m\otimes c_2\otimes c_1) (S\otimes 1\otimes S)[(
  \D(a)\otimes 1)\cdot\xi\cdot(1\otimes\D(b))]\no\\
& &~~~~~~~~~~=\e(a)\e(b)m(m\otimes 1)\cdot (1\otimes c_2\otimes c_1)
  (S\otimes 1\otimes S)\xi.
\eea
\end{Lemma}
The proof of this lemma is a straightforward computation, which we
omit.
\begin{Lemma}\label{L3} 
\bea
c_1^F&=&m\cdot(1\otimes c_1)(1\otimes S) F=\sum f_ic_1S(f^i),\no\\
c_2^F&=&m\cdot(1\otimes c_2)(S\otimes 1) F^{-1}=\sum 
    S(\bar{f}_i)c_2\bar{f}^i\label{twisted-c1-c2}
\eea
are invariant and pseudo-invariant, respectively, under the twisted
structure $\D_F(a)=F\D(a)F^{-1},~\forall a\in A$.
\end{Lemma}
\noindent{\bf Proof.} Write
\beq
\D_F(a)=\sum a^F_{(1)}\otimes a^F_{(2)}=F\sum a_{(1)}\otimes a_{(2)}
      F^{-1}.
\eeq
Then, $\forall a\in A$,
\bea
\sum S(a^F_{(1)}) c^F_2 a^F_{(2)}&=&\sum S(f_ia_{(1)}\bar{f}_j)
   S(\bar{f}_k)c_2\bar{f}^kf^ia_{(2)}\bar{f}^j
   (-1)^{[f^i]([a_{(1)}]+[\bar{f}_j])+[a_{(2)}][\bar{f}_j]}\no\\
&=&\sum S(\bar{f}_j)S(a_{(1)})S(\bar{f}_kf_i)c_2\bar{f}^k\bar{f}^i
   a_{(2)}\bar{f}^j(-1)^{[f_i][\bar{f}_k]+[a][\bar{f}_j]}\no\\
&=&\sum S(\bar{f}_j)S(a_{(1)})\e(F^{-1}F)c_2a_{(2)}\bar{f}^j
   (-1)^{[a][\bar{f}_j]}~~{\rm by}~(\ref{c-anti-inv})\no\\
&=&\e(a)\sum S(\bar{f}_j)c_2\bar{f}^j=\e(a)c_2^F.
\eea
Simlilarly, one can prove $\sum a^F_{(1)}c^F_1S(a^F_{(2)})=\e(a)c^F_1,~
\forall a\in A$.

We thus arrive at the following central elements induced by twisting
with $F$:
\bea
C^F_1&=&\sum S(X^F_\nu)\a_F Y^F_\nu c_1^FS(Z^F_\nu)=m(m\otimes 1)\cdot
   (1\otimes \a_F\otimes c_1^F)(S\otimes 1\otimes S)\Phi_F,\no\\
C^F_2&=&\sum S(X^F_\nu)c^F_2 Y^F_\nu\b_F S(Z^F_\nu)=m(m\otimes 1)\cdot
   (1\otimes c^F_2\otimes\b_F)(S\otimes 1\otimes S)\Phi_F\label{twisted-cbar}
\eea
which correspond to (\ref{C1-2}) and (\ref{C2-1}), respectively. Here
$\a_F$ and $\b_F$ are the twisted canonical elements given in 
(\ref{twisted-s-ab}).

\begin{Theorem}\label{univ-cbar}: The central elements
(\ref{twisted-cbar}) induced by twisting with $F$ coincide precisely
with the central elements $C_1,~C_2$ defined 
by (\ref{C1-2}) and (\ref{C2-1}), respectively. 
In other words, the central elements
$C_1$ and $C_2$ are invariant under twisting.
\end{Theorem}

To prove this theorem, we first notice
\begin{Lemma}\label{L4}: For any elements $\eta\in A\otimes A$ and
$\xi\in A\otimes A\otimes A$,
\bea
&(i)& m\cdot(1\otimes c_1^F)(1\otimes S)\eta=m\cdot(1\otimes c_1)
       (1\otimes S)(\eta F),\no\\
&(ii)& m\cdot(1\otimes c_2^F)(S\otimes 1)\eta=m\cdot(1\otimes c_2)
       (S\otimes 1)(F^{-1}\eta),\no\\
&(iii)& m\cdot(m\otimes 1)\cdot(1\otimes c_1^F\otimes c_2^F)
	(1\otimes S\otimes 1)\xi\no\\
& &~~~~=m\cdot(m\otimes 1)\cdot
	(1\otimes c_1\otimes c_2)(1\otimes S\otimes 1)
	[(1\otimes F^{-1})\cdot\xi\cdot(F\otimes 1)],\no\\
&(iv)& m\cdot(m\otimes 1)\cdot(1\otimes c_2^F\otimes
	c_1^F)(S\otimes 1\otimes S)\xi\no\\
& &~~~~=m\cdot(m\otimes 1)\cdot
	(1\otimes c_2\otimes c_1)(S\otimes 1\otimes S)
      [(F^{-1}\otimes 1)\cdot\xi\cdot(1\otimes F)].
\eea
\end{Lemma}
This lemma is proved by direct computation. Now with the help of
(\ref{twisted-phi}), lemmas \ref{L3},\ref{L4} and using the obvious
fact that $\b$, $\a$ are invariant and pseudo-invariant of $A$,
respectively, one can easily 
show that indeed $C^F_1=C_1$ and $C^F_2=C_2$.

For the quadratic-type invariants (\ref{quadra-c1}) and
(\ref{quadra-c2}), we have the central elements [c.f.
(\ref{quadra-cbars1})]
\bea
C_1&=&\sum S(X_\nu)\a Y_\nu\o_i\b S(\o^i)S(Z_\nu)
   =m(m\otimes 1)\cdot (S\otimes\a\otimes\b S)
   [\Phi(1\otimes\o)],\no\\
C_2&=&\sum S(X_\nu)S(\o_i)\a\o^i Y_\nu\b S(Z_\nu)
   =m(m\otimes 1)\cdot (S\otimes\a\otimes\b S)
   [(\o\otimes 1)\Phi].\label{quadra-cbars}
\eea
By (\ref{twisted-c1-c2}) and lemma \ref{L4}(i)(ii), one has
\bea
c^F_1&=&\sum f_jc_1S(f^j)=m\cdot(1\otimes \b S)(F\o)=m\cdot
   (1\otimes\b_F S)\o_F,\no\\
c^F_2&=&\sum S(\bar{f}_jc_2\bar{f}^j=m\cdot(S\otimes \a)(\o F^{-1})=m\cdot
   (S\otimes\a_F)\o_F,
\eea
where,
\beq
\o_F=F\o F^{-1}
\eeq
obviously commutes with the action of the twisted coproduct $\D_F$:
$\D_F(a)\o_F=\o_F\D_F(a),~~\forall a\in A$. In this notation,
we have central elements
\bea
C^F_1&=&m(m\otimes 1)\cdot
   (1\otimes \a_F\otimes \b_F)(S\otimes 1\otimes S)[\Phi_F(1\otimes
   \o_F)],\no\\
C^F_2&=&m(m\otimes 1)\cdot
   (1\otimes \a_F\otimes\b_F)(S\otimes 1\otimes S)[(\o_F\otimes 1)\Phi_F],
\eea
which, as a corollary of theorem \ref{univ-cbar}, 
reduce to $C_1$ and $C_2$ defined in
(\ref{quadra-cbars}),
respectively, independent of the twist applied.

In the case that $A$ is quasi-triangular with the universal R-matrix $\R$, 
where $\D(a)\R^T\R=\R^T\R\D(a),~\forall a\in A$, so we can take 
$\o=(\R^T\R)^m,~m\in {\bf Z}$. Then we obtain families of Casimir
invariants $C^m_1$ and $C^m_2,~m\in {\bf Z}$, which are invariant
under twisting.

\sect{Invariant Bilinear Forms and Invariant Forms}

Let $V, W$ be two (graded) $A$-modules, and $\ell(V,W)$ the space
of vector space maps (i.e. linear maps) from $V$ to $W$. We make
$\ell(V,W)$ into a (graded) $A$-module with the definition
\beq
(a\cdot f)(v)=\sum a_{(1)}f(S(a_{(2)})v)(-1)^{[f][a_{(2)}]},~~~~
\forall a\in A,~ v\in V,~ f\in\ell(V,W).
\eeq
We call $f$ invariant if
$a\cdot f=\e(a) f,~\forall a\in A$. Or equivalently, $\forall a\in A,~
v\in V$,
\beq
\sum a_{(1)}f(S(a_{(2)}v)(-1)^{[f][a_{(2)}]}=\e(a)f(v).\label{inv-f}
\eeq
In the case of normal Hopf (super)algebras, such invariants correspond
precisely to $A$-module homomorphisms, provided they are even. This is
not the case for quasi-Hopf (super)algebras. Nevertheless, there is
a close connection between such invariants and $A$-module
homomorphisms.
\begin{Proposition}\label{homo-f}: Suppose $f\in\ell(V,W)$ is even and
invariant. Set
\beq
\tilde{f}(v)=\sum S(X_\nu)\a Y_\nu f(S(Z_\nu)v),~~~~\forall v\in V.
\eeq
Then (i) $\tilde{f}\in \ell(V,W)$ is an $A$-module homomorphism, and
\bea
&(ii)& \b\tilde{f}(v)=f(v),~~~~\forall v\in V,\no\\
&(iii)& \tilde{f}(v)=\sum \bar{X}_\nu f(S(\bar{Y}_\nu)\a
  \bar{Z}_\nu v),~~~~\forall v\in V.
\eea
\end{Proposition}
\noindent{\bf Proof.}
Applying $m\cdot(1\otimes S)$ to lemma \ref{L1}(ii) and using
(\ref{inv-f}), one derives, 
\beq
\tilde{f}(S(a)v)=S(a)\tilde{f}(v),~~~~~\forall a\in A,~v\in V.
\eeq
Thus $\tilde{f}$ is an $A$-module homomorphism. This proves (i).
As for (ii), note
\beq
\tilde{f}(v)=m(m\otimes 1)\cdot (S\otimes\a\otimes f S)\cdot\Phi\cdot(1\otimes
   1\otimes v).
\eeq
Then by (\ref{quasi-bi2}),
\bea
1\otimes\tilde{f}(v)&=&\sum \bar{X}_\nu\bar{X}_\mu X^{(1)}_\s X_\rho\otimes
 S(Y_\rho)S(X^{(2)}_\s)S(\bar{Y}_\mu)S(\bar{Y}^{(1)}_\nu)\a
 \bar{Y}^{(2)}_\nu\bar{Z}_\mu Y_\s Z^{(1)}_\rho\no\\
& &\cdot f\lt(S(Z^{(2)}_\rho)S(Z_\s)S(\bar{Z}_\nu)v\rt) (-1)^y,
\eea
where
\beq
y=([X_\s]+[Z_\rho])([\bar{X}_\nu]+[\bar{X}_\mu])+[X^{(2)}_\s][Z_\rho]
  +[\bar{Z}_\nu][X_\s]+[\bar{X}_\mu][\bar{Y}_\nu]
  +[Z_\rho][X_\s].
\eeq
By (\ref{quasi-hopf1}) and (\ref{inv-f}),
\bea
1\otimes\tilde{f}(v)&=&\sum \bar{X}_\nu\bar{X}_\mu X^{(1)}_\s X_\rho\otimes
 S(Y_\rho)S(X^{(2)}_\s)S(\bar{Y}_\mu)\e(\bar{Y}_\nu)\a
 \bar{Z}_\mu Y_\s \e(Z_\rho)
 f(S(Z_\s)S(\bar{Z}_\nu)v) (-1)^y,\no\\
&=&(1\otimes m)\cdot(1\otimes S\otimes 1)
  \sum \lt(\bar{X}_\nu\bar{X}_\mu X^{(1)}_\s \otimes \e(\bar{Y}_\nu)
 \bar{Y}_\mu X^{(2)}_\s\otimes \a
 \bar{Z}_\mu Y_\s 
 f(S(Z_\s)S(\bar{Z}_\nu)v) \rt)\no\\
& & (X_\rho\otimes Y_\rho\otimes \e(Z_\rho))(-1)^{[X_\s]([\bar{Z}_\nu]
  +[\bar{X}_\mu])+[\bar{Z}_\nu][X_\s]+[\bar{Y}_\mu][X^{(2)}_\s]}\no\\
&=&  \sum \bar{X}_\nu\bar{X}_\mu X^{(1)}_\s \otimes \e(\bar{Y}_\nu)
  S(X^{(2)}_\s)S(\bar{Y}_\mu) \a
 \bar{Z}_\mu Y_\s 
 f\lt(S(\bar{Z}_\nu Z_\s)v\rt)\no\\
& & (-1)^{[\bar{Z}_\nu][Z_\s]
  +[X_\s][\bar{X}_\mu]}~~{\rm by}~(\ref{e(phi)=1})\no\\
&=&(1\otimes m)\cdot(1\otimes 1\otimes fS)\sum (
  \bar{X}_\nu\otimes \e(\bar{Y}_\nu)\otimes\bar{Z}_\nu)
  \lt(\bar{X}_\mu X^{(1)}_\s\rt.\no\\
& & \lt. \otimes S(X^{(2)}_\s)S(\bar{Y}_\mu)\a
  Z_\mu Y_\s\otimes Z_\s v\rt) (-1)^{[X_\s][\bar{X}_\mu]}\no\\
&=&\sum  \bar{X}_\mu X^{(1)}_\s\otimes S(X^{(2)}_\s)S(\bar{Y}_\mu)\a
  Z_\mu Y_\s f\lt(S(Z_\s) v\rt) (-1)^{[X_\s][\bar{X}_\mu]}~~
  {\rm by}~(\ref{quasi-bi4}).
\eea
Applying $m\cdot(1\otimes\b)$ gives
\bea
\b\tilde{f}(v)&=&
  \sum  \bar{X}_\mu X^{(1)}_\s\b S(X^{(2)}_\s)S(\bar{Y}_\mu)\a
  Z_\mu Y_\s f(S(Z_\s) v) (-1)^{[X_\s][\bar{X}_\mu]}\no\\
&=&\sum  \bar{X}_\mu \e(X_\s)\b S(\bar{Y}_\mu)\a
  Z_\mu Y_\s f(S(Z_\s) v) (-1)^{[X_\s][\bar{X}_\mu]}~~
  {\rm by}~(\ref{quasi-hopf2})\no\\
&=&\sum \e(X_\s)Y_\s f(S(Z_\s)v)=f(v)~~{\rm by}~(\ref{quasi-hopf3}),
  (\ref{e(phi)=1}),
\eea
which proves (ii). (iii) is a direct consequence of (ii) and (i).  

In the special case where $W={\bf C}$ is one-dimensional, 
we obtain the dual space
$V^*=\ell(V,{\bf C})$ which thus becomes a graded $A$-module with the
definition,
\beq
a\cdot f(v)=\sum\e(a_{(1)})f(S(a_{(2)})v)(-1)^{[f][a_{(2)}]}
  =(-1)^{[f][a]}f(S(a)v),~~~~\forall a\in A,~v\in V,~f\in V^*.
\eeq
We note that $f\in V^*$ is an $A$-invariant if and only if
\beq
\e(a)f(v)=a\cdot f(v)=(-1)^{[f][a]}f(S(a)v),~~~~\forall a\in A,
\eeq
or equivalently, since $\e(a)=0$ if $[a]=1$ and $\e(S^{-1}(a))=\e(a)$,
\beq
\e(a)f(v)=(S^{-1}\cdot f)(v)=f(av),~~~~~\forall a\in A.\label{inv-f*}
\eeq

A bilinear form $(~,~)$ on $V$ and $W$ is equivalent to an element
$\xi\in (V\otimes W)^*$ defined by
\beq
\xi(v\otimes w)=(v,w),~~~~~\forall v\in V,~w\in W.
\eeq
We say the form is invariant if $\xi$ is invariant. From
(\ref{inv-f*}) this is equivalent to
\beq
\e(a)\xi(v\otimes w)=\xi\lt(\D(a)(v\otimes w)\rt)=\sum\xi(a_{(1)}v\otimes
   a_{(2)}w)(-1)^{[v][a_{(2)}]},~~~~\forall a\in A.
\eeq
Thus a bilinear form is invariant iff
\beq
\sum (a_{(1)}v,a_{(2)}w)(-1)^{[v][a_{(2)}]}=\e(a)(v,w),~~~~\forall
  a\in A,~v\in V,~w\in W.
\eeq
In particular, a bilinear form $(~,~)$ on $A$ itself is called 
invariant iff
\beq
\sum({\rm Ad}a_{(1)}\cdot b,{\rm Ad}a_{(2)}\cdot c)(-1)^{[b][a_{(2)}]}
  =\e(a)(b,c),~~~~\forall a,b,c\in A.
\eeq

  Of particular interest are linear forms on $A$ which correspond to
elements $\xi$ of $A^*$. Such a linear form $\xi$ is called invariant
if it is an invariant element of $A^*$, i.e.
$a\cdot\xi=\e(a)\xi,~\forall a\in A$. Equivalently, $\xi\in A^*$ is
called an invariant linear form on $A^*$ if
\beq
\xi({\rm Ad}a\cdot b)=(S^{-1}(a)\cdot\xi)(b)=\e(a)\xi(b),~~~~\forall
   a,b\in A.\label{inv-xi}
\eeq
A linear form $\xi\in A^*$ is called pseudo-invariant if
\beq
\xi(\overline{\rm Ad}a\cdot b)=\e(a)\xi(b),~~~~\forall a,
  b\in A.\label{inv-xi-bar}
\eeq

Summarizing, $\xi\in A^*$ is called a pseudo-invariant linear form on $A$ if
\beq
\sum\xi(S(a_{(1)})ba_{(2)})(-1)^{[b][a_{(1)}]}=\e(a)\xi(b),~~~
   \forall a,b\in A,
\eeq
and an invariant linear form on $A$ if
\beq
\sum\xi(a_{(1)}bS(a_{(2)}))(-1)^{[b][a_{(2)}]}=\e(a)\xi(b),~~~~
   \forall a,b\in A.
\eeq
It is easily seen that
given any (graded) $A$-module $V$, the even invariants of $V^*=\ell
(V,{\bf C})$ correspond precisely with the $A$-module homomorphisms
$f\in V^*$. Thus the even invariant forms on $A$ correspond to
$A$-module homomorphisms $\xi\in A^*$, regarding $A$ as a module 
under the adjoint actions.

\sect{Casimir Invariants from Invariant Forms}

We now investigate the construction of central elements utilizing
invariant and pseudo-invariant linear forms on $A$. In the case
$A$ is quasi-triangular, we shall see how such central elements
may be constructed, corresponding to any finite dimensional
$A$-module, utilizing the universal R-matrix.

\begin{Proposition}\label{inv-linear-form}: Suppose
$\t=\sum a_i\otimes b_i\otimes c_i\in A^{\otimes 3}$ obeys
\beq
(1\otimes\D)\D(a)\cdot\t=\t\cdot(1\otimes\D)\D(a),~~~~\forall 
   a\in A.\label{t-d}
\eeq
If $\xi\in A^*$ is an even invariant form, then
\beq
C=(1\otimes\xi)(1\otimes m)(1\otimes 1\otimes\b S)\t=\sum a_i\xi
  (b_i\b S(c_i))
\eeq
is a central element. Similarly if $\bar{\t}=\sum\bar{a}_i\otimes
\bar{b}_i\otimes\bar{c}_i\in A^{\otimes 3}$ satisfies
\beq
\bar{\t}\cdot(\D\otimes 1)\D(a)=(\D\otimes 1)\D(a)\cdot\bar{\t},~~~~\forall
   a\in A,\label{tbar-d}
\eeq
and $\bar{\xi}\in A^*$ is an even pseudo-invariant form then
\beq
\bar{C}=(\bar{\xi}\otimes 1)(m\otimes 1)(S\otimes\a\otimes 1)\bar{\t}
  =\sum\bar{\xi}(S(\bar{a}_i)\a\bar{b}_i)\bar{c}_i
\eeq
is a central element.
\end{Proposition}
\noindent{\bf Proof.} 
Applying $(1\otimes m)(1\otimes 1\otimes \b S)$ to (\ref{t-d}),
one has
\bea
{\rm l.h.s.}&=&\sum a_{(1)}a_i\otimes {\rm Ad}a_{(2)}\cdot (b_i\b
  S(c_i))(-1)^{[a_i][a_{(2)}]},\no\\
{\rm r.h.s.}&=&\sum a_ia^R_{(1)}\otimes b_ia^R_{(2)}\b S(a^R_{(3)})
  S(c_i) (-1)^{[c_i]([a^R_{(2)}]+[a^R_{(3)}])+[a^R_{(1)}]([b_i]
  +[c_i])}\no\\
&=&\sum a_i a_{(1)}\e(a_{(2)})\otimes b_i\b S(c_i)(-1)^{[a]([b_i]
  +[c_i])}~~{\rm by}~(\ref{quasi-hopf2})\no\\
&=&\sum a_i a\otimes b_i\b S(c_i) (-1)^{[a]([b_i]+[c_i])}~~
  {\rm by}~(\ref{quasi-bi3}).
\eea
Applying $(1\otimes\xi)$ gives
\bea
{\rm l.h.s.}&=&\sum a_{(1)}a_i\otimes \xi\lt({\rm Ad}a_{(2)}\cdot (b_i\b
  S(c_i))\rt)(-1)^{[a_i][a_{(2)}]}\no\\
&=&\sum \lt(a_{(1)}a_i \xi\lt({\rm Ad}a_{(2)}\cdot (b_i\b
  S(c_i))\rt)\otimes 1\rt)(-1)^{[a_i][a_{(2)}]}\no\\
&=&\sum a_{(1)}a_i \e(a_{(2)})\xi(b_i\b
  S(c_i))\otimes 1=aC\otimes 1~~{\rm
  by}~(\ref{inv-xi}),(\ref{quasi-bi3}),\no\\
{\rm r.h.s.}&=&
   \sum a_i a\otimes \xi(b_i\b S(c_i)) (-1)^{[a]([b_i]+[c_i])}\no\\
&=& \sum a_i a \xi(b_i\b S(c_i))\otimes 1=Ca\otimes 1,
\eea
where in the second last equality
we have used the fact that $\xi$ is even, i.e. $\xi(a)=0$ if $[a]=1$.
This proves the first part of the proposition. The second part can be
proved in a similar way.       

It is easily shown that
\beq
\t=\Phi^{-1}(\o\otimes 1)\Phi,~~~~~~
\bar{\t}=\Phi(1\otimes \o)\Phi^{-1}\label{t-tbar}
\eeq
satisfy (\ref{t-d}), (\ref{tbar-d}), respectively. Thus as a 
corollary of proposition \ref{inv-linear-form} we have
the central elements,
\bea
C&=&\sum\xi(\bar{Y}_\nu\o^i Y_\mu\b S(\bar{Z}_\nu Z_\mu))\bar{X}_\nu
  \o_i X_\mu (-1)^{[X_\mu][\bar{Y}_\nu]+[Y_\mu][\bar{Z}_\nu]
  +[\o_i]([X_\mu]+[\bar{Y}_\nu])}\no\\
&=&(1\otimes\xi)(1\otimes m)(1\otimes 1\otimes\b S)[\Phi^{-1}(\o\otimes
  1)\Phi],\no\\
\bar{C}&=&\sum\bar{\xi}(S(X_\nu\bar{X}_\mu)\a Y_\nu\o_i \bar{Y}_\mu)Z_\nu
  \o^i \bar{Z}_\mu (-1)^{[Z_\nu][\bar{Z}_\mu]+[Y_\nu][\bar{X}_\mu]
  +[\o_i]([Z_\nu]+[\bar{Y}_\mu])}\no\\
&=&(\bar{\xi}\otimes 1)(m\otimes 1)(S\otimes \a\otimes 1)[\Phi(1\otimes
  \o)\Phi^{-1}].\label{inv-theta=o}
\eea

A quasi-Hopf (super)algebra is said to be of trace type if
there exists an invertible even
element $u\in A$ such that
\beq
S^2(a)=u au^{-1},~~~~\forall a\in A.\label{s2a=u}
\eeq
In the case $A$ is quasi-triangular with R-matrix as in (\ref{notation})
we have
\begin{Proposition}\label{u-operator}: The operator defined by
\beq
u=\sum S\lt(Y_\nu\b S(Z_\nu)\rt)S(e^i)\a e_i X_\nu (-1)^{[e_i]+
  [X_\nu]}\label{u}
\eeq
satisfies (\ref{s2a=u}). Moreover the inverse is given by
\beq
u^{-1}=S^2(u^{-1})=\sum S^{-1}(X_\nu)S^{-1}(\a\bar{e}^i)\bar{e}_i
  Y_\nu\b S(Z_\nu)(-1)^{[\bar{e}_i]}.\label{u-1}
\eeq
\end{Proposition}
\noindent {\bf Proof.} The non-super case was proved in \cite{Alt92}.
We here prove the super case. First observe
\bea
S^2(a)u&=&\sum S^2(a^L_{(3)})S(Y_\nu\b S(Z_\nu))S(e^i)S(a^L_{(1)})\a
  a^L_{(2)}e_iX_\nu\no\\
& &(-1)^{([a^L_{(1)}]+[a^L_{(2)}])([e_i]+[a^L_{(3)}]+[X_\nu])
  +[e_i]+[X_\nu]}~~{\rm by}~(\ref{quasi-hopf1}),(\ref{quasi-bi3})\no\\
&=&\sum S^2(a^L_{(3)})S(Y_\nu\b S(Z_\nu))m\cdot[(S\otimes \a)
  (a^L_{(1)}\otimes a^L_{(2)})(e^i\otimes e_i)]X_\nu\no\\
& &(-1)^{([a^L_{(1)}]+[a^L_{(2)}])([a^L_{(3)}]+[X_\nu])
  +[e_i]+[X_\nu]}\no\\
&=&\sum S^2(a^L_{(3)})S(Y_\nu\b S(Z_\nu))S(a^L_{(2)})S(e^i)\a
  e_ia^L_{(1)}X_\nu\no\\
& &(-1)^{([a^L_{(1)}]+[a^L_{(2)}])([a^L_{(3)}]+[X_\nu])
  +[e_i]+[X_\nu]+[a^L_{(1)}][a^L_{(2)}]}~~{\rm by}~(\ref{dr=rd})\no\\
&=&\sum S\lt(a^L_{(2)}Y_\nu\b S(Z_\nu)S(a^L_{(3)})\rt)S(e^i)\a e_i a^L_{(1)}
  X_\nu\no\\
& & (-1)^{[a^L_{(1)}]([a^L_{(2)}]+[a^L_{(3)}])
  +([a^L_{(1)}]+[a^L_{(3)}])[X_\nu] +[X_\nu]+[e_i]}.\label{1}
\eea
Applying $m\cdot T\cdot[\sum S(e^i)\a e_i (-1)^{[e_i]}\otimes S]$ to
lemma \ref{L1}(i), one has
\bea
ua&=&\sum S\lt(a^L_{(2)}Y_\nu\b S(Z_\nu)S(a^L_{(3)})\rt)S(e^i)\a e_i a^L_{(1)}
  X_\nu\no\\
& & (-1)^{[a^L_{(1)}]([a^L_{(2)}]+[a^L_{(3)}])
  +([a^L_{(1)}]+[a^L_{(3)}])[X_\nu] +[X_\nu]+[e_i]}\no\\
&=&S^2(a)u  ~~~{\rm by~(\ref{1})}.\label{s2au=ua}
\eea
It remains to show that $u$ is invertible. First we have
\begin{Lemma}\label{L5}:
\beq
S(\a)u=\sum S(e^i)\a e_i (-1)^{[e_i]}=m\cdot (S\otimes \a)\R^T.
\eeq
\end{Lemma}
\noindent{\bf Proof.} Note
\bea
u\otimes 1&=&[m\cdot(S\otimes\sum S(e^i)\a e_i(-1)^{[e_i]})(m\otimes 1)
  \otimes 1]\cdot (1\otimes\b S\otimes 1\otimes 1)\no\\
& &  (1\otimes T\otimes 1)
  (T\otimes 1\otimes 1) (\Phi\otimes 1)\no\\
&=&\sum S\lt(X^{(2)}_\nu Y_\mu\bar{X}_\s\bar{Y}^{(1)}_\rho\b
  S(\bar{Y}^{(2)}_\rho)S(Y_\nu Z^{(1)}_\mu\bar{Y}_\s)\rt)
  S(e^i)\a e_i X^{(1)}_\nu X_\mu\bar{X}_\rho\no\\
& &\otimes Z_\nu Z^{(2)}_\mu\bar{Z}_\s\bar{Z}_\rho (-1)^z~~{\rm by}~
  (\ref{quasi-bi2}),
\eea
where,
\bea
z&=&[e_i]+[X^{(1)}_\nu][X^{(2)}_\nu]+([X_\mu]+[X^{(1)}_\nu])
   ([Y_\mu]+[Y_\nu]+[\bar{Z}_\s]+[\bar{Y}_\rho]+[Z^{(1)}_\mu])\no\\
& &+[\bar{X}_\rho]([Z_\mu]+[Z_\nu]+[\bar{Z}_\s]+[\bar{Y}_\rho]+
   [Z^{(1)}_\mu])+[\bar{X}_\s]([Z_\mu]+[X_\nu])+[\bar{Y}_\s]
   ([Z_\nu]+[Z^{(2)}_\mu])\no\\
& &+[Z_\mu][X_\nu]+[Z_\nu][Z^{(1)}_\mu]+[\bar{Y}_\rho]
   ([Z_\mu]+[X_\nu]+[\bar{X}_\s]).
\eea
By (\ref{quasi-hopf2}), one has
\bea
u\otimes 1&=&
  \sum \lt[S\lt(X^{(2)}_\nu Y_\mu\bar{X}_\s\b
  S(Y_\nu Z^{(1)}_\mu\bar{Y}_\s)\rt)
  S(e^i)\a e_i X^{(1)}_\nu X_\mu
  \otimes Z_\nu Z^{(2)}_\mu\bar{Z}_\s\rt]\no\\
& & (\bar{X}_\rho\e(\bar{Y}_\rho)\otimes\bar{Z}_\rho) (-1)^{\bar{z}},
\eea
where
\bea
\bar{z}&=&[e_i]+[X^{(1)}_\nu][X^{(2)}_\nu]+([X_\mu]+[X^{(1)}_\nu])
  ([Y_\mu]+[Y_\nu]+[\bar{Z}_\s]+[Z^{(1)}_\mu])\no\\
& &+[\bar{X}_\s]([Z_\mu]+[X_\nu])+[\bar{Y}_\s]([Z_\nu]+Z^{(2)}_\mu])
  +[Z_\mu][X_\nu]+[Z_\nu][Z^{(1)}_\mu].
\eea
By (\ref{quasi-bi4}), one gets
\bea
u\otimes 1&=&
  \sum S\lt(X^{(2)}_\nu Y_\mu\bar{X}_\s\b
  S(Y_\nu Z^{(1)}_\mu\bar{Y}_\s)\rt)
  S(e^i)\a e_i X^{(1)}_\nu X_\mu
  \otimes Z_\nu Z^{(2)}_\mu\bar{Z}_\s
   (-1)^{\bar{z}}\no\\
&=&
  \sum S\lt(Y_\mu\bar{X}_\s\b
  S(Y_\nu Z^{(1)}_\mu\bar{Y}_\s)\rt)
  S(e^i)S(X^{(1)}_\nu)\a X^{(2)}_\nu e_iX_\mu
  \otimes Z_\nu Z^{(2)}_\mu\bar{Z}_\s\no\\
& &   (-1)^{\bar{z}+[X^{(2)}_\nu]([Y_\mu]+[Y_\nu]+[\bar{Z}_\s]
   +[Z^{(1)}_\mu])+[e_i]+[e_i][X_\nu]}~~{\rm by}~(\ref{dr=rd})\no\\
&=&\sum (\e(X_\nu)S^2(Y_\nu)\otimes X_\nu)
  \lt( S\lt(Y_\mu\bar{X}_\s\b
  S(Z^{(1)}_\mu\bar{Y}_\s)\rt)
  S(e^i)\a e_iX_\mu\rt.\no\\
& &\lt.  \otimes  Z^{(2)}_\mu\bar{Z}_\s\rt)
  (-1)^{[e_i]+[X_\mu]([Y_\mu]+[\bar{Z}_\s]
   +[Z^{(1)}_\mu])+[\bar{X}_\s][Z_\mu]+[\bar{Y}_\s][Z^{(2)}_\mu]}~~
   {\rm by}~(\ref{quasi-hopf1})\no\\
&=&\sum 
   S\lt(Y_\mu\bar{X}_\s\b S(\bar{Y}_\s)
  S(Z^{(1)}_\mu)\rt)
  S(e^i)\a e_iX_\mu
  \otimes Z^{(2)}_\mu\bar{Z}_\s\no\\
& &  (-1)^{[e_i]+[X_\mu]([Y_\mu]+[\bar{Z}_\s]
   +[Z^{(1)}_\mu])+[\bar{X}_\s][Z_\mu]+[\bar{Y}_\s][Z_\mu]}~~
   {\rm by}~(\ref{e(phi)=1}).
\eea
Applying $m\cdot(1\otimes S(\a))\cdot T\cdot (1\otimes S)$ gives 
\bea
S(\a)u&=&\sum S(\bar{Z}_\s)S(S(Z^{(1)}_\mu)\a Z^{(2)}_\mu)
   S(Y_\mu\bar{X}_\s\b S(\bar{Y}_\s)) S(e^i)\a e_iX_\mu\no\\
& &(-1)^{[Z_\mu]+[\bar{Z}_\s]+[e_i]+[X_\mu]([Y_\mu]+[\bar{Z}_\s])}\no\\
&=&\sum S(\bar{Z}_\s)\e(Z_\mu)S(\a)
   S(Y_\mu\bar{X}_\s\b S(\bar{Y}_\s)) S(e^i)\a e_iX_\mu\no\\
& &(-1)^{[\bar{Z}_\s]+[e_i]+[X_\mu]([Y_\mu]+[\bar{Z}_\s])}~~{\rm by}~
  (\ref{quasi-hopf1})\no\\
&=&\sum S(\bar{Z}_\s) S(\a) S(\bar{X}_\s\b S(\bar{Y}_\s))
  (-1)^{[\bar{Z}_\s]}\no\\
& &m\cdot T\cdot (\sum S(e^i)\a e_i (-1)^{[e_i]}\otimes S)
  (X_\mu\otimes Y_\mu\e(Z_\mu))\no\\
&=&\sum S(\bar{Z}_\s) S(\a) S(\bar{X}_\s\b S(\bar{Y}_\s))
   S(e^i)\a e_i (-1)^{[e_i]+[\bar{Z}_\s]}~~{\rm
   by}~(\ref{e(phi)=1})\no\\
&=&\sum S(\bar{X}_\s\b S(\bar{Y}_\s)\a\bar{Z}_\s)
   S(e^i)\a e_i (-1)^{[e_i]}\no\\
&=&\sum S(e^i)\a e_i (-1)^{[e_i]}~~{\rm by}~(\ref{quasi-hopf3}),
\eea
thus proving lemma \ref{L5}.
\begin{Lemma}\label{L6}: 
\beq
u\sum S^{-1}(\a\bar{e}^i)\bar{e}_i (-1)^{[e_i]}=\a.
\eeq
\end{Lemma}
This lemma is easily proved with the help of (\ref{s2au=ua}) and lemma
\ref{L5}.

Now we are in a position to prove proposition \ref{u-operator}:
\bea
1&=&\sum S(X_\nu)\a Y_\nu\b S(Z_\nu)~~{\rm by}~(\ref{quasi-hopf4})\no\\
&=&\sum S(X_\nu)u S^{-1}(\a\bar{e}^i)\bar{e}_iY_\nu\b
   S(Z_\nu)(-1)^{[\bar{e}_i]}~~{\rm by~lemma}~\ref{L6}\no\\
&=&u\sum S^{-1}(X_\nu)S^{-1}(\a\bar{e}^i)\bar{e}_iY_\nu
  \b S(Z_\nu)(-1)^{[\bar{e}_i]}~~{\rm by}~(\ref{s2au=ua})\no\\
&=&S^2\lt(\sum S^{-1}(X_\nu)S^{-1}(\a\bar{e}^i)\bar{e}_iY_\nu
  \b S(Z_\nu)(-1)^{[\bar{e}_i]}\rt) u~~{\rm by}~(\ref{s2au=ua}).
\eea
It follows that $u$ is invertible, with $u^{-1}$ given by (\ref{u-1}).
This completes our proof for proposition \ref{u-operator}.    

\begin{Corollary}: If $A$ is a quasi-tringular quasi-Hopf (super)algebra,
then $A$ is of trace-type.
In particular $S^2(u)=u$, and $uS(u)=S(u)u$ is a central element.
\end{Corollary}

Below we assume that $A$ is a quasi-Hopf (super)algebra of trace type.
Let $V$ be finite-dimensional (graded) $A$-module. Then

\begin{Proposition}\label{trace-form}: $\xi\in A^*$ defined by
\beq
\xi(a)={\rm Str}_V(uS^{-1}(\a)a),~~~~\forall a\in A
\eeq
determines an invariant linear form, and $\bar{\xi}\in A^*$
defined by
\beq
\bar{\xi}(a)={\rm Str}_V(u^{-1}S(\b)a),~~~~\forall a\in A
\eeq
determines a pseudo-invariant linear form.
\end{Proposition}
\noindent{\bf Proof.} By means of (\ref{s2a=u}) and the (super)trace property
\beq
{\rm Str}_V(ab)=(-1)^{[a][b]}{\rm Str}_V (ba),~~~~\forall a,b
   \in A\label{super-trace}
\eeq 
one has, $\forall a,b\in A$,
\bea
\xi({\rm Ad}a\cdot b)&=&\sum {\rm Str}_V\lt(u S^{-1}(\a)a_{(1)}b
  S(a_{(2)})\rt)(-1)^{[b][a_{(2)}]}\no\\
&=&\sum {\rm Str}_V\lt(S(a_{(2)}u S^{-1}(\a)a_{(1)}b
  \rt)(-1)^{[a_{(1)}][a_{(2)}]}~~{\rm by}~(\ref{super-trace})\no\\
&=&\sum {\rm Str}_V\lt(uS^{-1}(a_{(2)} S^{-1}(\a)a_{(1)}b
  \rt)(-1)^{[a_{(1)}][a_{(2)}]}~~{\rm by}~(\ref{s2a=u})\no\\
&=&\sum {\rm Str}_V\lt(uS^{-1}(S(a_{(1)})\a a_{(2)})b\rt)\no\\
&=&\e(a)\cdot {\rm Str}_V(uS^{-1}(\a)b)~~{\rm by}~(\ref{quasi-hopf1})\no\\
&=&\e(a)\xi(b).
\eea
Thus we have proved the first part of the proposition.
The second part of the proposition is proved in a similar fashion.   

It immediately follows from propositions \ref{inv-linear-form} and 
\ref{trace-form} that one has

\begin{Proposition}\label{trace-inv}: Let $\pi$ be the representation
afforded by the finite-dimensional (graded) $A$-module $V$. Suppose
$\t=\sum a_i\otimes B_i\otimes c_i\in A\otimes {\rm End}V\otimes A$ obeys
\beq
(1\otimes\pi\otimes 1)(1\otimes\D)\D(a)\cdot\t=\t\cdot
 (1\otimes\pi\otimes 1)(1\otimes\D)\D(a),~~~~\forall a\in A,
\eeq
then
\beq
C=\sum {\rm Str}_V\lt(uS^{-1}(\a)B_i\b S(c_i)\rt)a_i
\eeq
is a central element. Similarly if $\bar{\t}=\sum\bar{a}_i\otimes
\bar{B}_i\otimes\bar{c}_i\in A\otimes {\rm End}V\otimes A$ satisfies
\beq
\bar{\t}\cdot(1\otimes\pi\otimes 1)(\D\otimes 1)\D(a)=
  (1\otimes\pi\otimes 1)(\D\otimes 1)\D(a)\cdot\bar{\t},~~~~\forall
   a\in A.
\eeq
Then
\beq
\bar{C}=\sum {\rm Str}_V\lt(u^{-1}S(\b)S(\bar{a}_i)\a\bar{B}_i\rt)\bar{c}_i
\eeq
is a central element.
\end{Proposition}

\begin{Corollary}: Suppose $\o=\sum \o_i\otimes \Omega^i\in A\otimes{\rm
End}V$ satisfies
\beq
(1\otimes\pi)\D(a)\cdot\o=\o\cdot(1\otimes\pi)\D(a),~~~~\forall a\in A.
\eeq
Then the first equation of (\ref{inv-theta=o}) implies that
\bea
C&=&=\sum{\rm Str}_V\lt(uS^{-1}(\a)\bar{Y}_\nu\Omega^iY_\mu\b S(\bar{Z}_\nu
  Z_\mu)\rt)\bar{X}_\nu\o_i X_\mu\no\\
& &(-1)^{[X_\mu][\bar{X}_\nu]+[Y_\mu][\bar{Y}_\nu]+[\o_i]
  ([X_\mu]+\bar{Y}_\nu])}\label{trace-o-inv}
\eea
is a central element. Similarly if $\bar{\o}=\sum\bar{\Omega}_i\otimes
\bar{\o}^i\in {\rm End}V\otimes A$ satisfies
\beq
(\pi\otimes 1)\D(a)\cdot\bar{\o}=\bar{\o}\cdot(\pi\otimes 1)\D(a),~~~~
  \forall a\in A.
\eeq
Then the second equation of (\ref{inv-theta=o}) means that
\bea
\bar{C}&=&\sum{\rm Str}_V\lt(u^{-1}S(\b)S(X_\nu\bar{X}_\mu)\a Y_\nu
  \bar{\Omega}_i\bar{Y}_\mu\rt)Z_\nu\bar{\o}^i\bar{Z}_\mu\no\\
& &(-1)^{[Z_\nu][\bar{Z}_\mu]+[Y_\nu][\bar{X}_\mu]+
  [\o^i]([Z_\nu]+[\bar{Y}_\mu])}\label{trace-obar-inv}
\eea
is a central element.
\end{Corollary}

\begin{Corollary}
In the case that $A$ is quasi-triangular, one takes $\o=(\R^T\R)^m,~
m\in {\bf Z}$. Then we obtain the following families of Casimir
invariants associated with $\R^T\R$ and its powers:
\bea
C_m&=&(1\otimes{\rm Str}_V)(1\otimes m)(1\otimes uS^{-1}(\a)\otimes\b
   S)\cdot\Phi^{-1}((\R^T\R)^m\otimes 1)\Phi,\no\\
\bar{C}_m&=&({\rm Str}_V\otimes 1)(m\otimes 1)(u^{-1}S(\b)S\otimes\a
   \otimes 1)\cdot\Phi(1\otimes (\R^T\R)^m)\Phi^{-1}.\label{cm}
\eea
\end{Corollary}
\noindent{\bf Remark:} The above invariants are natural generalizations
of those obtained in \cite{Lin92,Zha91}
to which they reduce in the case of normal
Hopf (super)algebras (for which $\Phi=1\otimes 1\otimes 1$).

\sect{Twisting Invariance of Central Elements $C$ and $\bar{C}$}

In this section we show that the trace-type central elements $C$ and
$\bar{C}$ are invariant under twisting. Associated with $F$,
we have the twisted co-associator $\Phi_F$ and in the
quasi-triangular case, the twisted R-matrix $\R_F$.
We write,
\bea
\Phi_F&=&\sum X^F_\nu\otimes Y^F_\nu\otimes Z^F_\nu,\no\\
\Phi_F^{-1}&=&\sum \bar{X}^F_\nu\otimes \bar{Y}^F_\nu\otimes 
   \bar{Z}^F_\nu,\no\\
\R_F&=&\sum e^F_t\otimes e^t_F.
\eea
\begin{Lemma}\label{L7}
\beq
\b=\sum\bar{f}_i\b_F S(\bar{f}^i),~~~~
\a=\sum S(f_i)\a_F f^i.
\eeq
\end{Lemma}
This lemma is proved by direct computation by means of
(\ref{twisted-s-ab}).

Associated with a twist $F$ on a quasi-triangular quasi-Hopf
(super)algebra, we have the $u$-operator in terms of
the twisted structure, denoted $u_F$:
\beq
u_F=\sum S\lt(Y^F_\nu\b_F S(Z^F_\nu)\rt)S(e_F^t)\a_F 
  e^F_t X^F_\nu (-1)^{[e^F_t]
  +[X^F_\nu]}\label{uF}
\eeq

\begin{Theorem}\label{inv-u}: The $u$-operator, given explicitly in
proposition \ref{u-operator}, is invariant under twisting.
\end{Theorem}
\noindent{\bf Proof.} We compute $u_F$. By (\ref{twisted-phi}), one has
\bea
u_F&=&\sum S\lt(f^if_j^{(2)}Y_\nu\bar{f}^k_{(1)}\bar{f}_\l\b_F
  S(\bar{f}^l)S(\bar{f}^k_{(2)}S(f^j Z_\nu)\rt)\no\\
& &S(e^t_F)\a_F e^F_tf_if^{(1)}_jX_\nu\bar{f}_k (-1)^{[e^F_t]+r},
\eea
where
\bea
r&=& [f_i]+([f^{(1)}_j]+X_\nu]+[\bar{f}_k])
  ([f^i]+[f^{(2)}_j]+[f^j]+[\bar{f}^k]+[X_\nu])
  +[f^j]([\bar{f}_k]+[Z_\nu])\no\\
& &+[f^k]([f^j]+[Z_\nu])+[f^{(2)}_j]([\bar{f}_k]+[X_\nu])
  +[\bar{f}_k][X_\nu].
\eea
By lemma \ref{L7},
\bea
u_F&=&\sum S\lt(f^if_j^{(2)}Y_\nu\bar{f}^k_{(1)}\b
  S(\bar{f}^k_{(2)}S(f^j Z_\nu)\rt)\no\\
& &S(e^t_F)\a_F e^F_tf_if^{(1)}_jX_\nu\bar{f}_k (-1)^{[e^F_t]+r},\no\\
&=&\sum S\lt(f^if_j^{(2)}Y_\nu\e(\bar{f}^k)\b
  S(f^j Z_\nu)\rt)
  S(e^t_F)\a_F e^F_tf_if^{(1)}_jX_\nu\bar{f}_k\no\\
& & (-1)^{[e^F_t]+[f_i]+[X_\nu]+[f^j][Y_\nu]+[X_\nu][f^i]
  +[f^{(1)}_j]([f^{(2)}_j]+[f^i]+[f^j]+[X_\nu])}~~{\rm
  by}~(\ref{quasi-hopf2})\no\\
&=&\sum S\lt(f_j^{(2)}Y_\nu\b
  S(f^j Z_\nu)\rt)
  S(e^t_Ff^i)\a_F e^F_tf_if^{(1)}_jX_\nu\no\\
& & (-1)^{[e^F_t]+[f_i]+[X_\nu]+[f^i][e^F_t]+[f^j][Y_\nu]
  +[f^{(1)}_j]([f^{(2)}_j]+[f^j]+[X_\nu])}~~{\rm
  by}~(\ref{e(f)=1})\no\\
&=&\sum S\lt(f_j^{(2)}Y_\nu\b
  S(f^j Z_\nu)\rt)
  S(e^t)S(f^i)\a_F f^ie_tf^{(1)}_jX_\nu\no\\
& & (-1)^{[e_t]+[X_\nu]+[f^j][Y_\nu]
  +[f^{(1)}_j]([f^{(2)}_j]+[f^j]+[X_\nu])}~~{\rm
  by}~(\ref{dr=rd})\no\\
&=&\sum S\lt(f_j^{(2)}Y_\nu\b
  S(f^j Z_\nu)\rt)
  S(e^t)\a e_tf^{(1)}_jX_\nu\no\\
& & (-1)^{[e_t]+[X_\nu]+[f^j][Y_\nu]
  +[f^{(1)}_j]([f^{(2)}_j]+[f^j]+[X_\nu])}~~{\rm by~lemma}~\ref{L7}\no\\
&=&\sum S\lt(Y_\nu\b
  S(f^j Z_\nu)\rt)
  S(e^tf^{(2)}_j)\a e_tf^{(1)}_jX_\nu\no\\
& & (-1)^{[e_t]+[X_\nu]+[f^{(2)}_j][e_t]
  +[f^{(1)}_j][f^{(2)}_j]+[f_j]+[f_j][Z_\nu]}\no\\
&=&\sum S\lt(Y_\nu\b
  S(f^j Z_\nu)\rt)
  S(e^t)S(f^{(1)}_j)\a f^{(2)}_je_tX_\nu\no\\
& & (-1)^{[e_t]+[X_\nu]+[f_j][e_t]
  +[f_j]+[f_j][Z_\nu]}~~{\rm by}~(\ref{dr=rd})\no\\
&=&\sum S\lt(Y_\nu\b
  S(f^j Z_\nu)\rt)
  S(e^t)\a e_t\e(f_j)X_\nu
 (-1)^{[e_t]+[X_\nu]}~~{\rm by}~(\ref{quasi-hopf1})\no\\
&=&\sum S\lt(Y_\nu\b
  S(Z_\nu)\rt)
  S(e^t)\a e_tX_\nu
 (-1)^{[e_t]+[X_\nu]}~~{\rm by}~(\ref{e(f)=1})\no\\
&=&u.
\eea
Thus we end up with the same $u$-operator, independently of the
twist applied. 
\begin{Corollary}:
\beq
S(u)S(\b)=\sum e_i\b S(e^i)=m\cdot(\b\otimes S)R.
\eeq
\end{Corollary}
\noindent{\bf Proof:} We apply theorem 4 and lemma \ref{L5} to the
special case where $F$ is the Drinfeld twist \cite{Dri90} $F_D$.
In \cite{Gou98}, we proved
\beq
S(\b)=\a_{F_D},~~~~(S\otimes S)\R=\R_{F_D}.\label{a'r'}
\eeq
Then from lemma \ref{L5} and theorem 4,
\beq
S(\a_{F_D})u=m\cdot (S\otimes \a_{F_D})\R^T_{F_D}
\eeq
which gives rise to, on using (\ref{a'r'}),
\beq
S^2(\b)u=\sum S^2(e^i)S(\b)S(e_i) (-1)^{[e_i]}=S\sum e_i\b S(e^i).
\eeq
Namely,
\beq
\sum e_i\b S(e^i)=S^{-1}(u)S(\b)=S(u)S(\b),
\eeq
where we have used $S^2(u)=u\cdot u\cdot u^{-1}=u$.

The following result follows as a special case of proposition 
\ref{trace-form} applied to the twisted quasi-Hopf (super)algebra
structure.
\begin{Lemma}\label{L8}: $\xi_F\in A^*$ defined by
\beq
\xi_F(a)={\rm Str}_V(uS^{-1}(\a_F)a),~~~~\forall a\in A
\eeq
determines a linear form invariant under the twisted quasi-Hopf (super)
algebra structure. Similarly $\bar{\xi}\in A^*$ defined by
\beq
\bar{\xi}_F(a)={\rm Str}_V(u^{-1}S(\b_F)a),~~~~\forall a\in A
\eeq
determines a pseudo-invariant linear form under the twisted structure.
\end{Lemma}

Following proposition \ref{inv-linear-form}, if $\t\in A^{\otimes 3}$
satisfies (\ref{t-d}) and $\bar{\t}\in A^{\otimes 3}$ satisfies 
(\ref{tbar-d}), then we have trace type invariants
\bea
C&=&(1\otimes{\rm Str})(1\otimes m)(1\otimes uS^{-1}(\a)\otimes \b
   S)\t,\no\\
\bar{C}&=&({\rm Str}\otimes 1)(m\otimes 1)(u^{-1}S(\b) S\otimes\a\otimes
   1)\bar{\t}.\label{trace-c-cbar}
\eea

\begin{Lemma}\label{L9}: Suppose $\t\in A^{\otimes 3}$ satisfies 
(\ref{t-d}). Then
\beq
\t_F\equiv (1\otimes F)(1\otimes\D)F\cdot\t\cdot(1\otimes\D)F^{-1}
   (1\otimes F^{-1})\label{tF}
\eeq
also satisfies (\ref{t-d}) for the twisted structure; viz
\beq
(1\otimes\D_F)\D_F(a)\cdot\t_F=\t_F\cdot(1\otimes\D_F)\D_F(a),~~~~
   \forall a\in A.\label{d-tF}
\eeq
Similarly if $\bar{\t}\in A^{\otimes 3}$ satisfies 
(\ref{tbar-d}). Then
\beq
\bar{\t}_F\equiv (F\otimes 1)(\D\otimes 1)F\cdot\bar{\t}\cdot(\D\otimes
  1)F^{-1} (F^{-1}\otimes 1)\label{tbarF}
\eeq
also satisfies (\ref{tbar-d}) for the twisted structure; viz
\beq
(\D_F\otimes 1)\D_F(a)\cdot\bar{\t}_F=\bar{\t}_F\cdot(\D_F\otimes 1)\D_F(a),~~~~
   \forall a\in A.\label{d-tbarF}
\eeq
\end{Lemma}
\noindent{\bf Proof.} Applying $(1\otimes F)(1\otimes\D)F$ to the left
and $(1\otimes\D)F^{-1}(1\otimes F^{-1})$ to the right of (\ref{t-d})
gives (\ref{d-tF}). Similarly, applying $(F\otimes 1)(\D\otimes 1)F$ to
the left and $(\D\otimes 1)F^{-1}(F^{-1}\otimes 1)$ to the right
of (\ref{tbar-d}), one gets (\ref{d-tbarF}).

We thus arrive at the following central elements obtained by twisting 
those of (\ref{trace-c-cbar}) with $F$:
\bea
C_F&=&(1\otimes{\rm Str})(1\otimes m)(1\otimes uS^{-1}(\a_F)\otimes \b_F
   S)\t_F,\no\\
\bar{C}_F&=&({\rm Str}\otimes 1)(m\otimes 1)(u^{-1}S(\b_F) S\otimes\a_F\otimes
   1)\bar{\t}_F.\label{trace-c-cbar-F}
\eea
We shall show that these invariants coincide precisely with those of
(\ref{trace-c-cbar}). Namely, 

\begin{Theorem}\label{inv-trace}: The trace type central elements 
(\ref{trace-c-cbar}) are invariant under twisting. 
\end{Theorem}

To prove this theorem, we first state
\begin{Lemma}\label{L10}: $\forall a\in A,~\xi\in A^{\otimes 3}$,
we have
\bea
&(i)& (1\otimes{\rm Str}_V)(1\otimes m)(1\otimes uS^{-1}(\a)\otimes
  \b S)\cdot\xi(1\otimes\D(a))\no\\
& & ~~~~~~~~=(1\otimes{\rm Str}_V)(1\otimes m)(1\otimes uS^{-1}(\a)\otimes
  \b S)\cdot(1\otimes\D(a))\xi\no\\
& & ~~~~~~~~=\e(a)\;(1\otimes{\rm Str}_V)(1\otimes m)(1\otimes uS^{-1}(\a)
  \otimes\b S)\cdot\xi\no\\
&(ii)& ({\rm Str}_V\otimes 1)(m\otimes 1)(u^{-1}S(\b)S\otimes\a\otimes
  1)\cdot(\D(a)\otimes 1)\xi\no\\
& &~~~~~~~~= ({\rm Str}_V\otimes 1)(m\otimes 1)(u^{-1}S(\b)S\otimes\a\otimes
  1)\cdot\xi(\D(a)\otimes 1)\no\\
& &~~~~~~~~=\e(a)\; ({\rm Str}_V\otimes 1)(m\otimes 1)(u^{-1}S(\b)S\otimes\a
  \otimes 1)\cdot\xi\no\\
&(iii)& (1\otimes{\rm Str}_V)(1\otimes m)(1\otimes uS^{-1}(\a_F)
  \otimes\b_F S)\cdot\xi\no\\
& &~~~~~~~~=(1\otimes{\rm Str}_V)(1\otimes m)(1\otimes uS^{-1}(\a)
  \otimes\b S)[(1\otimes F^{-1})\cdot\xi\cdot(1\otimes F)]\no\\
&(iv)& ({\rm Str}_V\otimes 1)(m\otimes 1)(u^{-1}S(\b_F)S\otimes
  \a_F\otimes 1)\cdot\xi\no\\
& & ({\rm Str}_V\otimes 1)(m\otimes 1)(u^{-1}S(\b)S\otimes
  \a\otimes 1)[(F^{-1}\otimes 1)\cdot\xi\cdot(F\otimes 1)].
\eea
\end{Lemma}
\noindent{\bf Proof.} This lemma is proved  by direct computations using
(\ref{super-trace}), (\ref{s2a=u}), (\ref{quasi-hopf1}) and 
(\ref{quasi-hopf2}). For demonstration, we show the details for proving
some of the relations. Write $\xi=\sum x_i\otimes y_i\otimes z_i\in
A^{\otimes 3}$. Then,
\bea
&&(1\otimes{\rm Str}_V)(1\otimes m)(1\otimes uS^{-1}(\a)\otimes
  \b S)\cdot(1\otimes\D(a))\xi\no\\
&&~~~~~~=\sum(1\otimes{\rm Str}_V)\lt(x_i\otimes uS^{-1}(\a)a_{(1)}
   y_i\b S(z_i)S(a_{(2)})\rt) (-1)^{[x_i]([a]+[a_{(2)}])}\no\\
&&~~~~~~=\sum(1\otimes{\rm Str}_V)\lt(x_i\otimes S(a_{(2)})uS^{-1}(\a)a_{(1)}
   y_i\b S(z_i))\rt)\no\\
&&~~~~~~~~~~   (-1)^{[x_i][a]+[a_{(1)}][a_{(2)}]}~~{\rm by}~
   (\ref{super-trace})\no\\
&&~~~~~~=\sum(1\otimes{\rm Str}_V)\lt(x_i\otimes uS^{-1}(S(a_{(1)}\a
   a_{(2)})y_i\b S(z_i))\rt) (-1)^{[x_i][a]}~~{\rm by}~
   (\ref{s2a=u})\no\\
&&~~~~~~=\e(a)\;(1\otimes{\rm Str}_V)(1\otimes m)(1\otimes uS^{-1}(\a)
  \otimes\b S)\cdot\xi~~{\rm by}~(\ref{quasi-hopf1}).
\eea
Other relations in (i) and (ii) are proved similarly. We now prove (iii):
\bea
&& (1\otimes{\rm Str}_V)(1\otimes m)(1\otimes uS^{-1}(\a_F)
  \otimes\b_F S)\cdot\xi\no\\
&&~~~~~~=\sum(1\otimes{\rm Str}_V)\lt(x_i \otimes uS^{-1}(\bar{f}^j)S^{-1}(\a)
    \bar{f}_jy_if_k\b S(z_if^k)\rt)\no\\
&&~~~~~~~~~~    (-1)^{[\bar{f}_j]+[z_i][f_k]}~~{\rm by}~
    (\ref{twisted-s-ab})\no\\
&&~~~~~~=\sum(1\otimes{\rm Str}_V)\lt(x_i \otimes S(\bar{f}^j)uS^{-1}(\a)
    \bar{f}_jy_if_k\b S(z_if^k)\rt)\no\\
&&~~~~~~~~~~    (-1)^{[\bar{f}_j]+[z_i][f_k]}~~{\rm by}~
    (\ref{s2a=u})\no\\
&&~~~~~~=\sum(1\otimes{\rm Str}_V)\lt(x_i \otimes uS^{-1}(\a)
    \bar{f}_jy_if_k\b S(\bar{f}^jz_if^k)\rt)\no\\
&&~~~~~~~~~~    (-1)^{[\bar{f}_j]([y_i]+[f_k])+[z_i][f_k]}~~{\rm by}~
    (\ref{super-trace})\no\\
&&~~~~~~=(1\otimes{\rm Str}_V)(1\otimes m)(1\otimes uS^{-1}(\a)
  \otimes\b S)[(1\otimes F^{-1})\cdot\xi\cdot(1\otimes F)].
\eea
(iv) can be proved in a similar fashion.

We are now in a position to prove theorem \ref{inv-trace}. From
(\ref{trace-c-cbar-F}), one has, by lemma \ref{L10}(iii) and (\ref{tF}),
\bea
C_F&=&(1\otimes{\rm Str})(1\otimes m)(1\otimes uS^{-1}(\a)\otimes \b
   S)[(1\otimes\D)F\cdot\t\cdot(1\otimes\D) F^{-1}]\no\\
&=&(1\otimes{\rm Str})(1\otimes m)(1\otimes uS^{-1}(\a)\otimes \b
   S)[(1\otimes\e)F\cdot\t\cdot(1\otimes\e) F^{-1}]~~{\rm by~
   lemma}~\ref{L10}(i)\no\\
&=&(1\otimes{\rm Str})(1\otimes m)(1\otimes uS^{-1}(\a)\otimes \b
   S)\t~~{\rm by}~(\ref{e(f)=1})=C.
\eea
Similarly, one can show $\bar{C}_F=\bar{C}$. This completes the proof
of theorem \ref{inv-trace}.

In the quasi-triangular case it is worth noting that when 
$\t,~\bar{\t}$ have the special form of (\ref{t-tbar}) with 
$\o=(\R^T\R)^m,~\in m\in{\bf Z}$, then their twisted analogues are
given by
\beq
\t_F=\Phi^{-1}_F(\o_F\otimes 1)\Phi_F,~~~~
\bar{\t}_F=\Phi_F(1\otimes \o_F)\Phi_F^{-1},~~~~\o_F=(\R^T_F\R_F)^m,
   \label{theta-F}
\eeq
which agree precisely with the prescription of lemma \ref{L9}. It
follows, as a special case of theorem 5, that the central elements
of (\ref{cm}) are invariant under twisting.

\vskip.3in
\noindent {\bf Acknowledgements.}
The financial support from Australian Research 
Council through a Queen Elizabeth II Fellowship Grant for Y.-Z.Z is
gratefully acknowledged. P.S.I has been financially supported by
an Australian Postgraduate Award.

\newpage

\end{document}